\documentclass{amsart}

\usepackage{fullpage,amsmath,amsfonts}
\usepackage[all]{xy}
\usepackage{amssymb,amsthm}
\usepackage{hyperref}

\newcommand{\Set}{\mathbf{Set}}

\newcommand{\calC}{\ensuremath{\mathcal{C}}}
\newcommand{\calD}{\ensuremath{\mathcal{D}}}
\newcommand{\calE}{\ensuremath{\mathcal{E}}}

\newcommand{\calI}{\ensuremath{\mathcal{I}}}
\newcommand{\calJ}{\ensuremath{\mathcal{J}}}
\newcommand{\calL}{\ensuremath{\mathcal{L}}}
\newcommand{\calM}{\ensuremath{\mathcal{M}}}
\newcommand{\calS}{\ensuremath{\mathcal{S}}}

\newcommand{\opCat}[1]{\ensuremath{#1^\mathrm{op}}}

\newcommand{\Nat}{\ensuremath{\mathbb{N}}}

\newcommand{\Psh}[1]{\widehat{#1}}

\newcommand{\rmC}{\mathrm{C}}  %% A comonad

\theoremstyle{plain}
\newtheorem{theorem}{Theorem}[section]
\newtheorem{corollary}[theorem]{Corollary}
\newtheorem{proposition}[theorem]{Proposition}
\newtheorem{lemma}[theorem]{Lemma}

\theoremstyle{definition}
\newtheorem{definition}[theorem]{Definition}
\newtheorem{example}[theorem]{Example}
\newtheorem{remark}[theorem]{Remark}

\begin{document}

%opening
\title{The successive dimension, without elegance}

\author{Mat\'\i as Menni}
\address{Conicet and  Departamento de Matem\'atica, Universidad Nacional de La Plata, Argentina}
\email{matias.menni@gmail.com}
\urladdr{sites.google.com/site/matiasmenni/} % Delete if not wanted.
\thanks{This work was supported by CONICET (Argentina), PIP 11220200100912CO}

\begin{abstract}
Experience shows that the poset of levels (or dimensions) of the topos of presheaves on some elegant Reedy categories  may be equipped with a monotone increasing `successor' function which, as the case of simplicial sets shows, is different from Lawvere's Aufhebung in general.
We prove that a similar result holds for the topos of presheaves on a small category with split-epi/mono factorizations; a typical feature of categories that are  Reedy elegant, or  skeletal, or  graphic (von Neumann-)regular, but more general.
In fact, we show that the more general `successor' may be described as a function on the poset of   full subcategories of the site that are closed under subobjects.
\end{abstract}

\maketitle

%\tableofcontents

\section{Introduction}

Hurewicz and Wallman start  \cite{DimensionTheory}  by celebrating Poincar\'e who, writing in a philosophical journal and ``concerned only with putting forth and intuitive concept of dimension and not an exact mathematical formulation [...]  had, however, penetrated very deep, in stressing the inductive nature of the geometric meaning of dimension and the possibility of disconnecting a space by subsets of lower dimension''. 
Immediately after, they sketch their precise definition (attributed to Menger and Urysohn): the empty set has dimension $-1$; the dimension of a (separable metric) space is the least integer $n$ for which every point has arbitrarily small neighborhoods whose boundaries have dimension less than $n$.

For some `tame' spaces the dimension has an alternative characterization. For example consider the following inductively defined Heyting formulas over variables ${ \{X_d \mid d \in \Nat \}}$:
\begin{equation*}
    B_d =
    \begin{cases}
      X_0 \vee \neg X_0  , & \text{if}\ d=0 \\
      ( X_d \vee (X_d \Rightarrow B_{d-1} ), & \text{if}\ d \geq 1.
    \end{cases}
\end{equation*}
 These are  the {\em Bounded Depth} formulas. Among other results, \cite[Theorem~{4.1}]{MarraEtAl2018} shows that a non-empty polyhedron $P$ has dimension $d$ (in the sense of  \cite{DimensionTheory}) if and only if, in the Heyting algebra of open sub-polyhedra of $P$, ${B_d}$ holds and ${B_{c}}$ fails for every ${0 \leq c < d}$.

Some alternative definitions  of `dimension' are discussed in \cite[Section~{I.5}]{DimensionTheory} but, in all cases, dimension is an integer or a real number.
The Dimension Theory for toposes outlined in \cite[Section~{II}]{Lawvere91} is different; we recall below some of the basic details.

A {\em level} of  a topos ${\calE}$ is a string  ${l_! \dashv l^* \dashv l_* : \calL \rightarrow \calE}$ with fully faithful ${l_! , l_*  }$.
So $\calL$ is a topos and ${l : \calL \rightarrow \calE}$ is an essential subtopos with direct image $l_*$.
Quoting from \cite{Lawvere91}: ``The basic idea is simply to identify dimensions with levels and then try to determine what the general dimensions are in particular examples. More precisely, a space may be said to have (less than or equal to) the dimension grasped by a given
level if it belongs to the negative (left adjoint inclusion) incarnation of that level.'' 
For this reason, the leftmost adjoint $l_!$  may be called the {\em ($l$-)skeletal inclusion} and, for each $X$ in $\calE$, the counit ${l_! (l^* X) \rightarrow X}$ will be called the {\em ($l$-)skeleton} of $X$. So $X$ is {\em $l$-skeletal} if its $l$-skeleton is an isomorphism.

Levels of a  given topos may be partially ordered as subtoposes. 
So each topos, but especially each topos of spaces, determines its poset of dimensions.
Notice that there is nothing `inductive' about this.
For example,   the poset of levels of the topos ${\Psh{\Delta}}$ of simplicial sets  may be identified with ${\{ -\infty, 0, 1, \ldots, n, \ldots, \infty \}}$ but, in general, the poset of levels of a topos need not be totally ordered. See \cite{KellyLawvere89, Lawvere2002a, MMlevelEpsilon} and references therein. 

Lawvere observes that, in many examples, the poset of levels/dimensions may be equipped with a monotone  increasing function called {\em Aufhebung} which assigns, to each level, the smallest level ``qualitatively higher'' than it. More precisely, for  levels $l$, $m$ of a topos $\calE$,  $m$ is {\em way-above}  $l$ if ${l \leq m}$ (that is, as subtoposes) and also ${l_! : \calL \rightarrow \calE}$ factors through ${m_* : \calM \rightarrow \calE}$.
The Aufhebung of $l$ is, when it exists, the least level that is way-above $l$. It is not easy to calculate with the tools we have today. Fairly general techniques for this task were developed for toposes of presheaves on graphic monoids \cite{Lawvere2002a}; the case of ball complexes is solved in \cite{Roy2021}; the case of simplicial sets and similar examples is treated in \cite{KennetEtAl2011};  the case of the classifier of non-trivial Boolean algebras remains an open problem \cite{LawvereOpenProblems}; the case of the `gros' Zariski topos and other toposes in Algebraic Geometry seem out of reach at present.
(It is relevant to observe that, in the case of $\Psh{\Delta}$, the Aufhebung  is not the obvious successor function: for ${n \geq 2}$, the Aufhebung of level~$n$ is ${2 n - 1}$.)

As we have already stressed, no induction over the natural numbers is needed to define the poset of dimensions of a topos but, of course, we may iterate the Aufhebung. Each topos has a degenerate level (called ${-\infty}$) which is the only essential subtopos whose domain is terminal as a category.
Its Aufhebung is called {\em level~0} and the Aufhebung of level 0 is called {\em level~1}.
For instance, level~0 of ${\Psh{\Delta}}$  is ${\Set = \Psh{\Delta_0} \rightarrow \Psh{\Delta}}$ and level~1 is ${\Psh{\Delta_1} \rightarrow \Psh{\Delta}}$. 

Although the framework is quite different from  \cite{DimensionTheory}, the geometric intuition remains.
If the $0$-skeletal inclusion  has an additional left adjoint (connected-components), then level~1  coincides with  ``the smallest dimension such that the set of components of an arbitrary space is the same as the set of components of the skeleton at that dimension of the space''. (See \cite[Proposition in p.~19]{Roy2021} or \cite[Corollary~{2.5}]{MMlevelEpsilon} for a proof.)  
Hence, as in \cite{DimensionTheory},  `dimension~1'  is related with connectedness of the space.

Let ${l : \calL \rightarrow \calE}$ be a level. For simplicity let us assume that the skeleton ${\ell_X : l_! (l^* X) \rightarrow X}$ is monic for each $X$ in $\calE$. A subobject ${u : U \rightarrow X}$ is said to have {\em $\ell$-skeletal boundary} if ${id_X \leq u \vee (u \Rightarrow \ell_X)}$ in the Heyting algebra  of subobjects of $X$. 
(Notice that, if that  algebra was also co-Heyting, as in any presheaf topos, then $u$ having  $\ell$-skeletal boundary would be equivalent to  ${\partial u \leq \ell_X}$, where ${\partial u}$ denotes the usual `boundary' operation in a co-Heyting algebra. Compare with the Bounded Depth formulas defined above.)

With $\calE$ and $l$ as above,  an object $X$ in $\calE$ is said to have {\em $l$-skeletal boundaries} if every subobject of $X$ has skeletal boundary. The objects with $l$-skeletal boundaries seem to be  closely related with the objects that are skeletal with respect to a qualitatively higher level (but not necessarily its Aufhebung). The relation is not fully understood but, for example,  for certain presheaf toposes, 1-skeletal objects coincide with quotients of 0-separated objects with 0-skeletal  boundaries  \cite{Menni2019a}. The sites to which these results apply are typical elegant Reedy categories as defined in \cite{BergnerRezk2013}. Technically, `elegance' guarantees that  skeleta  are monic. More general results, as the ones in this paper, must deal with that issue differently; `without elegance', so it speak, in the sense that we do not assume  that absolute pushouts of split epimorphisms exist in the site.

The skeletal counit of any level may be factored as an epic followed by a monic and this determines  an idempotent comonad with monic counit that we call the {\em principal} comonad induced by the level.
In Section~\ref{SecCoend} we characterize the coalgebras for principal comonads  induced by  levels in presheaf toposes. 
This characterization has independent interest and  is fundamental for the calculations in the rest of the paper.

Experience with toposes of presheaves on elegant Reedy categories suggests to concentrate on the `special'  levels such that the counit of the principal comonad is mono-cartesian.  In Section~\ref{SecMonoCartesianIdeals} we show that, if $\calD$ is a small category with split-epic/mono factorizations, then the idempotent ideals of $\calD$ that induce such special levels are in bijective correspondence with the full subcategories of $\calD$ that are closed under subobjects.
(The existence of split-epi/mono factorizations is a typical feature of sites in Algebraic Topology and Homotopy Theory.
For example, the elegant Reedy categories already mentioned, the skeletal categories in the sense of \cite[Definition~{8.1.1}]{Cisinski2006}, but also the graphic (von Neumann-)regular categories used to study the Aufhebung in \cite[Proposition~2]{Lawvere2002a}.)

In Section~\ref{SecMapsOnTopOfAmonoCoreflection} we show that every idempotent comonad with monic  and mono-cartesian counit on a topos determines an idempotent  ideal of the topos. For the topos  of presheaves on a site $\calD$ with split-epi/mono factorizations,  that ideal restricts to an idempotent ideal (of $\calD$) which may be thought of as the `successor' of the original comonad. Restricting to principal comonads induced by idempotent ideals in the site determines a `successor' function on the poset of special levels which, in view of the above, may be identified with a function  on the poset of full subcategories of $\calD$ that are closed under subobjects.
A description of this `successor', purely in terms of $\calD$, is given in Section~\ref{SecMapsOnTopOfPrincipalComonad}.
We discuss several examples in Section~\ref{SecExamples}.

\section{Principal comonads of levels in presheaf toposes}
\label{SecCoend}

Each level of a topos determines an idempotent comonad with monic counit that we call the {\em principal} comonad associated to the level.
In some cases, such as simplicial sets and other toposes of presheaves on elegant Reedy categories, this comonad coincides with the skeletal comonad of the level; but this is not so in general. The purpose of the present section is to give a description of the principal comonad for a level in an arbitrary presheaf topos. We first quickly recall the necessary background from \cite{KellyLawvere89} and \cite{Menni2021}.

Let $\calE$ be a topos. 
A (Lawvere-Tierney) topology  $j$ in $\calE$  is called {\em principal} if each object has a least $j$-dense subobject.
For example, if  ${\calE_j \rightarrow \calE}$ is a level then $j$ is principal.
Indeed, the image  of the counit ${j_! (j^* X) \rightarrow X}$ is the least $j$-dense subobject of $X$.

Let $j$ be a principal topology in $\calE$ and, for each $X$ in $\calE$, let ${\gamma_X : \rmC X \rightarrow X}$ be the least $j$-dense subobject.
For future reference we emphasize the following.

\begin{definition}\label{DefMinimal}
The object  $X$ is $j$-minimal if $\gamma_X$ is an isomorphism.
\end{definition}

The full subcategory of $j$-minimal objects is coreflective and ${\gamma_X : \rmC X \rightarrow X}$ is universal from this subcategory to $X$.
In other words, ${(\rmC, \gamma)}$ is  an idempotent comonad whose coalgebras are the $j$-minimal objects.
This comonad will be called the {\em principal comonad} determined by the (principal) topology $j$.
See \cite[Section~3]{Menni2021} for further details.

We next recall the relation between levels of a presheaf topos and idempotent ideals in the site.
A {\em right ideal} of a category $\calC$ is a family ${ ( \calI(-, C) \mid  C \in \calC)}$ where each ${\calI(-, C)}$ is a sieve on $C$.
For ${f : B \rightarrow C}$ in $\calC$, we may write ${f \in \calI}$ instead of ${f \in \calI(B, C)}$.
A {\em (two-sided) ideal} is a right  ideal $\calI$ such that if ${g f}$ is defined and ${f \in \calI}$ then ${g f \in \calI}$.
An ideal $\calI$ is {\em idempotent} if, for every ${f \in \calI}$ there are ${g, h \in \calI}$ such that ${f = g h}$.

We stress that, in the definition above, $\calC$ need not be small.
On the other hand, for small $\calC$, \cite[Theorem~{4.4}]{KellyLawvere89} proves that there is monotone bijection between levels of the presheaf topos ${\Psh{\calC}}$ and idempotent ideals of $\calC$. If $\calI$ is an idempotent  ideal of $\calC$ then the corresponding (essential) subtopos of ${\Psh{\calC}}$ is determined by the Grothendieck topology on $\calC$ such that a sieve covers an object $C$ if it contains the sieve ${\calI(-,C)}$.

Fix an idempotent ideal $\calI$ in a small category $\calD$.

Recall that for a Grothendieck topology on a small category, the sheaf reflection of a presheaf   may be built by applying the +-construction twice
\cite[Section~{III.5}]{maclane2}.
For topologies determined by  ideals, the +-construction may be simplified.

\begin{lemma}\label{LemPlusConstructionForIdeals} 
For the Grothendieck topology induced by the   ideal $\calI$, the following holds: for any presheaf $P$ in ${\Psh{\calD}}$ and $D$  in $\calD$, ${P^{+} D}$ may be identified with ${\Psh{\calD}(\calI(-,D), P)}$ or, alternatively, with the set of families
\[ (x_f \in P C \mid f   \in \calI(C, D), C \in \calD ) \]
 that are {\em matching}: ${x_f \cdot g = x_{f g}}$ for any $f$ in $\calI$ and any $g$ in $\calD$ post-composable with $f$.
\end{lemma}
\begin{proof}
In general, ${P^+ D}$ is the set of equivalence classes of matching families (indexed by the maps in a covering sieve of $D$).
Two  families being equivalent if they coincide in a common (covering) refinement of the corresponding sieves.
In our present case, ${\calI D}$  is the least covering refinement of any cover of $D$,
so any equivalence class of matching families (over $D$) is uniquely determined by their refinement to that least covering.
\end{proof}

A  family as in  Lemma~\ref{LemPlusConstructionForIdeals} may be written as ${(x_f \mid f \in \calI D) \in P^{+} D}$ and then
 \[ (x_f \mid f \in \calI D) \cdot h = (x_{h k} \mid k \in \calI C ) \in P^{+} C \]
holds for every ${h : C \rightarrow D}$ in $\calD$.

There is a canonical  ${P \rightarrow P^+}$ which, in the context of Lemma~\ref{LemPlusConstructionForIdeals}, may be described, for each $D$, as the function ${P D \rightarrow P^{+} D}$ that sends ${x \in P D}$ to the matching family ${(x\cdot f \mid f \in \calI  D) \in P^+ D}$.
As mentioned above, the composite ${P \rightarrow P^+ \rightarrow  (P^+)^+}$ 
is universal from $P$ to  sheaves. The sheaf reflection $(P^+)^+$ is sometimes denoted by  ${\mathbf{a} P}$.

In order to prove the main result of the section we will need a concrete description of the skeletal inclusion determined by the idempotent ideal $\calI$.
To give such a description we rely on wedges and (co)ends as discussed in \cite[Chapter~{IX}]{maclane}.

Let ${T : \calD \times \opCat{\calD} \times \calD \rightarrow \Psh{\calD}}$ be the functor 
\[ (D, C, B) \mapsto  \calI(C, D) \times \calI(-, B) \]
 and let us consider ${T_D = T(D, -, -) : \opCat{\calD} \times \calD \rightarrow \Psh{\calD}}$

 \begin{lemma}\label{LemCompatibleFamiliesAsWedges} For any $D$ in $\calD$ and any $P$ in $\Psh{\calD}$, the set 
 \[ \Psh{\calD}\left( \int^C T_D(C,C), P \right) \]
 is isomorphic with ${(\mathbf{a}P) D}$. Moreover, this isomorphism is natural in $D$.
 \end{lemma}
 \begin{proof}
 We first show that the set of  wedges from ${T_D}$ to $P$ is isomorphic to ${(\mathbf{a} P) D}$.
 Since ${\mathbf{a} P = (P^+)^+}$, it follows that  ${(\mathbf{a} P) D}$ may be identified with the set of matching families
\[ ( \mathbf{x}_f \in P^+ C \mid f \in \calI(C, D) , C\in \calD ) \]
so that, if for ${f : C \rightarrow D}$ in $\calI$ we let ${ \mathbf{x}_f = ( x_{f, g} \mid g  \in \calI(B, C), B \in \calD ) }$, then the matching property says that 
\[ (x_{f, h k} \mid  k \in \calI B) = (x_{f, g} \mid g \in \calI C) \cdot h  = (\mathbf{x}_f) \cdot h = \mathbf{x}_{f h} =  (x_{f h, k} \mid k \in \calI B) \]
for every  ${h : B  \rightarrow C}$ in $\calD$.
 
On the other hand,  a wedge from ${T_D}$ to $P$  is a family  ${(\omega_B \mid B \in \calD)}$ of maps 
 \[ \xymatrix{ T_D (C, C) = \calI(C, D) \times \calI(-, C) \ar[r]^-{\omega_C} & P} \]
 in $\Psh{\calD}$ such that the diagram
 \[\xymatrix{
 \calI(C, D) \times \calI(-, B) \ar[d]_-{\calI(h, C) \times \calI(-,B)} \ar[rrr]^-{\calI(C, D) \times \calI(-,h)} &&&  \calI(C, D) \times \calI(-, C) \ar[d]^-{\omega_C} \\
  \calI(B, D) \times \calI(-, B) \ar[rrr]_-{\omega_B} &&& P
 }\]
 commutes for every ${h : B \rightarrow C}$ in $\calD$. 
 
 If, for ${f  \in  \calI(C, D)}$ and ${g \in \calI( B, C) }$, we let ${\omega_{f, g} = \omega_C (f, g) \in P B}$ then,  the wedge $\omega$ may be identified with a family
 \[ ( (\omega_{f, g} \in P B \mid  g \in \calI(B, C) , B\in \calD ) \mid f  \in \calI(C, D) , C \in \calD) \]
satisfying two naturality conditions.  The first one says that the family
\[ (\omega_{f, g} \in P B \mid  g \in \calI(B, C) , B\in \calD ) \]
is matching for every  ${f : C \rightarrow D}$. The second one (i.e. the wedge property) says that
\[   \omega_{f h, k} = \omega_{f, h k} \]
for every ${h : B \rightarrow C}$ in $\calD$ and every ${k \in \calI B}$.
So  the set of  wedges from ${T_D}$ to $P$ is indeed isomorphic to ${(\mathbf{a} P) D}$.
Hence, the iso in the statement follows from the universal property of coends.
Naturality follows from the `Parameter Theorem'  \cite[Theorem~{IX.7.2}]{maclane}.
 \end{proof}

The proof of \cite[Theorem~{4.4}]{KellyLawvere89}, showing that essential subtoposes of a presheaf topos are in bijective correspondence with ideals in the site, rests on Freyd's SAFT to produce the left adjoint of the inverse image of a subtopos determined  by an ideal. Here we give an alternative description that will relevant for the rest of the paper.

\begin{proposition}[Coend formula for skeleta]\label{PropSkeletaAsCoend}
The endofunctor on ${\Psh{\calD}}$ defined by  
\[Q \mapsto Q_! = \int^{C, D} Q D \times \calI(C, D) \times \calI(-, C) \]
 extends to a left adjoint to sheafification.
\end{proposition}
\begin{proof}
Let $Q$ be a presheaf on $\calD$ and let 
${S_Q : \opCat{\calD} \times \calD \times \opCat{\calD} \times \calD \rightarrow \Psh{\calD}}$ be  defined by
\[ S_Q (D, C, B, A) =  Q D \times T_C (B, A)   = Q D \times T(C, B, A) = Q D \times \calI(B, C) \times \calI(-, A) \]
so that  ${ Q_!  = \int^{C, D} S_Q(D, D, C, C)}$. Then, for any presheaf  $P$  calculate
\begin{align*}
\Psh{\calD}(Q_!, P) &= \Psh{\calD}\left( \int^{C, D} S_Q (D, D, C , C) , P \right)  \\
 &\cong \int_{C, D} \Psh{\calD}(Q D \times T_D(C,C), P) & \text{Continuity} \\
 &\cong \int_{C,D} \int_B \Set( Q D \times (T_D(C, C))B, P B)   & \text{end formula for Nat} \\
 &\cong  \int_{C,D} \int_B \Set( Q D, \Set ((T_D(C, C))B, P B)) & \text{Exponential transposition} \\
 &\cong \int_D \Set \left(Q D, \int_C \int_B \Set ((T_D(C, C))B, P B)\right) &  \text{Fubini and Continuity} \\
 &\cong \int_D \Set\left(Q D, \int_C \Psh{\calD}(T_D(C,C), P) \right)  & \text{end formula for Nat} \\
 & \cong \int_D \Set\left(Q D, \Psh{\calD}\left(\int^C T_D(C,C), P\right) \right)  & \text{Continuity} \\
 &\cong \int_D \Set( Q D, (\mathbf{a} P) D)  & \text{Lemma~\ref{LemCompatibleFamiliesAsWedges}} \\
 &\cong\Psh{\calD}(Q, \mathbf{a} P)  & \text{end formula for Nat}.
\end{align*}
\end{proof}

Let $P$ be a presheaf on $\calD$. For any $D$ in $\calD$, we say  that  ${x \in P D}$ is  {\em $\calI$-generated} if there exists an ${f : D \rightarrow E}$ in $\calI$ and a ${y \in P E}$ such that ${x = y \cdot f}$. The presheaf  $P$  is said to be {\em $\calI$-generated} if,  for every $D$ in $\calD$ and ${x \in P D}$,  $x$ is  $\calI$-generated in the sense of the previous sentence.

\begin{example}[The case of representable presheaves]\label{ExIgeneratedCore} If we let ${P = \calD(-,C)}$ for some $C$ in $\calD$ then ${x \in P D}$ is $\calI$-generated if and only if ${x \in \calI}$.
Indeed, by definition, $x$ is $\calI$-generated if and only if there is an ${f : D \rightarrow E}$ in $\calI$ and a ${y: E \rightarrow C}$ in $\calD$ such that ${x = y \cdot f}$.
The left  ideal property implies that ${x \in \calI}$. Conversely, if ${x \in \calI}$ then idempotency implies the existence of ${y, f \in \calI}$ such that ${x = y \cdot f}$. 
Hence, $P$ is $\calI$-generated if and only if ${id_C \in \calI}$.
\end{example}

The next result compares $\calI$-generated presheaves  with the minimal ones in the sense of Definition~\ref{DefMinimal}.

\begin{theorem}\label{ThmCoendFormula} Let $\calD$ be a small category and let $l$ be the level of $\Psh{\calD}$ induced by an idempotent ideal $\calI$ in $\calD$.
A presheaf on $\calD$ is $l$-minimal if and only if it is $\calI$-generated. 
\end{theorem}
\begin{proof}
We described  the unit ${P \rightarrow \mathbf{a} P}$  after Lemma~\ref{LemPlusConstructionForIdeals}.
Taking it as an element  in ${\Psh{\calD}(P, \mathbf{a} P)}$ at the bottom of the calculation in the proof of Proposition~\ref{PropSkeletaAsCoend}, and chasing it `up', we may conclude that its transposition ${P_! \rightarrow P}$ is determined by the wedge
\[ \xymatrix{
P D \times \calI(C, D) \times \calI(-, C) \rightarrow P
}\]
that sends ${x \in P D}$ and ${f : C \rightarrow D}$, ${g : B \rightarrow C}$ in $\calI$ to ${x \cdot f \cdot g = x \cdot (f g) \in P B}$.
It follows that the image of  ${P_! B  \rightarrow P B}$ consists of the elements of the form ${x \cdot (f g) \in P B}$ for some $x$, $f$ and $g$ as above but, since $\calI$ is idempotent, these are the elements of the form ${x \cdot h}$ for some ${h : B \rightarrow D}$ in $\calI$ and ${x \in P D}$.
\end{proof}

In other words, the principal comonad on ${\Psh{\calD}}$ associated to  (the level determined by) the idempotent ideal $\calI$ sends a presheaf to the subpresheaf  of $\calI$-generated elements. In particular, by Example~\ref{ExIgeneratedCore}, the counit at stage representable by $C$ in $\calD$ is ${\calI(-,C) \rightarrow \calD(-,C)}$.

\section{Mono-cartesian ideals}
\label{SecMonoCartesianIdeals}

A natural transformation is {\em mono-cartesian} if, for every monic in the relevant category, the corresponding naturality square is a pullback.
We are concerned with the special case of a monic natural transformation that is the counit of an idempotent comonad.
If the underlying category has pullbacks then such a monic counit is mono-cartesian if and only if the  subcategory of coalgebras if closed under  subobjects \cite[Lemma~{2.6}]{Menni2021}.

We are interested in mono-cartesian principal comonads on presheaf toposes.
Again, let $\calI$ be an idempotent ideal of a small category $\calD$.

\begin{lemma}\label{LemMonoCartesian} 
If the idempotent ideal  ${\calI}$ in $\calD$ is such that:
\begin{itemize}
\item for every  map  ${f : D \rightarrow E}$ in $\calI$  there exists a factorization ${f = g b}$ with ${b : D \rightarrow B}$ and a  map ${s : B \rightarrow D}$ such that ${s b \in \calI}$ and ${b s b = b}$,
\end{itemize}
then the induced principal comonad is mono-cartesian.
\end{lemma}
\begin{proof}
We rely on the description of the principal comonad given by Theorem~\ref{ThmCoendFormula}.
Let ${m : X \rightarrow Y}$ be monic in $\Psh{\calD}$. To prove that the naturality  square below
\[\xymatrix{
\rmC X \ar[d]_-{\gamma_X} \ar[rr]^-{\rmC m} && \rmC  Y \ar[d]^-{\gamma_Y} \\
X \ar[rr]_-{m} &&  Y
}\]
is a pullback, let ${x \in X D}$ and assume that ${m_D x }$ is $\calI$-generated. That is, there is a ${f : D \rightarrow E}$ in $\calI$  and a ${z \in Y E}$ such that ${m_D x = z \cdot f}$. Let $g$, $b$ and $s$ as in the statement. Then ${m_D x = (z\cdot g) \cdot b}$ and also
\[  m_D (x \cdot s \cdot b) = (m_D x) \cdot s  \cdot b = (z\cdot g) \cdot  b \cdot s \cdot b = (z\cdot g) \cdot b = m_D x  \]
and infer, using that $m$ is monic, that ${(x \cdot s) \cdot b = x}$.
So $x$ is $\calI$-generated. This completes the proof that $\gamma$ is mono-cartesian.
\end{proof}

Notice that Lemma~\ref{LemMonoCartesian} is applicable if ${f \in \calI}$ factors as ${f = g b}$ with $b$ a split epic in $\calI$.

Naturally, an idempotent ideal is called {\em mono-cartesian} if the corresponding principal comonad has mono-cartesian counit.

\begin{proposition}\label{PropMonoCartesianIdeals} 
If every map in the small $\calD$ factors as a split epic followed by a monomorphism then, the mono-cartesian idempotent ideals in $\calD$ are in bijective correspondence with the full  subcategories of $\calD$ that are closed under subobjects.
\end{proposition}
\begin{proof}
First observe that, by the right-ideal property, if a map in $\calI$ factors as split epic followed by monic then the monic part is also in $\calI$.
Also, if the monic ${m : D \rightarrow E}$ is in $\calI$ then the mono-cartesian property implies that ${id_D}$ is in $\calI$.
(Just use the mono-cartesian property with the monic ${\calD(-,m) : \calD(-,D) \rightarrow \calD(-, E)}$.)

Let $\calI$ be a mono-cartesian idempotent  ideal in $\calD$. Let ${\calC \rightarrow \calD}$ be the full subcategory consisting of those objects  whose corresponding identity is in $\calI$,  and let ${\calJ}$ be the associated idempotent  ideal in $\calD$, that is, the ideal of maps that factor through some object in $\calC$. Certainly, ${\calJ \subseteq \calI}$ because if ${f = g h}$ with ${g : C\rightarrow D}$,
$C$ in $\calC$ then, as ${id_C \in \calI}$, ${f = g (id_C ) h}$ must also be in the two-sided ideal $\calI$.

To prove that ${\calI \subseteq \calJ}$ let ${ f : D\rightarrow F}$ be in $\calI$. 
By hypothesis we can factor $f$ as ${f = m e}$ with
${m : E \rightarrow F}$  monic and ${e : D\rightarrow E}$ split-epic. 
As observed above, it follows that $m$ is in $\calI$ and, by mono-cartesianness,  that $id_E$ is in $\calI$.
So $f$ is in $\calJ$.

Now let $C$ be in the  subcategory ${\calC \rightarrow \calD}$   and let ${m : D \rightarrow C}$ be monic in $\calD$.
Then ${ id_C \in \calI}$ and, by the right-ideal property, ${m \in \calI}$. Since we are assuming that the ideal is mono-cartesian, ${id_D \in\calI}$, so $D$ is in $\calC$.

Conversely, assume that $\calI$ is the idempotent ideal determined by a full subcategory ${\calC \rightarrow \calD}$ closed under subobjects.
 Let ${f : D \rightarrow E}$ in $\calI$.
That is, ${f = u v}$  with ${u : C \rightarrow E}$ and $C$ in $\calC$.
Also, by hypothesis, ${v : D \rightarrow C}$ factors as a split epic ${b : D \rightarrow B}$ followed by a monic ${m : B \rightarrow C}$.
Then $B$ is in $\calC$ by hypothesis, so ${id_B \in \calI}$ and, therefore, ${b \in \calI}$.
Taking ${g = u m}$ we may apply Lemma~\ref{LemMonoCartesian} to conclude that $\calI$ is mono-cartesian.
\end{proof}

As a corollary we may deduce the following strengthening of  \cite[Proposition~{2.5}]{Menni2019a}.

\begin{corollary}\label{PropEnoughSplittingsImplyMonoCartesianness} 
If the small  $\calD$ is such that  every map factors as a split epic followed by a split monic then the  levels of ${\Psh{\calD}}$ are in bijective correspondence with the  full  subcategories of $\calD$ that are closed under subobjects. If, moreover, every object of $\calD$ has a finite set of subobjects then every subtopos is a level.
\end{corollary}

\section{The ideal of maps on top of a mono-coreflection}
\label{SecMapsOnTopOfAmonoCoreflection}

Let $H$ co-Heyting algebra with top element ${\top}$ and  ${(-)/u \dashv u \vee (-)}$ for each  ${u \in H}$.
Then ${\partial u = (\top/u) \wedge u}$ is  sometimes called the {\em  boundary} of $u$ (e.g. in \cite{Lawvere91a}).
Notice that, if $H$ is also a Heyting algebra then  ${\partial u \leq v \in H}$ if and only if ${\top \leq u \vee (u \Rightarrow v)}$.
Notice also that the second inequality makes sense in any Heyting category.

The co-Heyting boundary can be generalized as follows. For any ${f \in H}$ we may define ${\partial_f u = f/u \wedge u}$.
Again, in a bi-Heyting algebra,  ${\partial_f u \leq v}$ if and only if ${f \leq u \vee  (u \Rightarrow v)}$.
So, for any ${f, u, v}$ in a Heyting algebra $H$ we may write ${\partial_f u \leq v}$ if ${f \leq u \vee  (u \Rightarrow v)}$.
To avoid a possible confusion we stress that  ``${\partial_f u}$'' does not make sense in this generality;  only the expresions of the form ${\partial_f u \leq v}$ do.

We assume from now on an underlying Heyting category. The implication in the Heyting algebra of subobjects of a fixed object will be denoted by $\Rightarrow$.

\begin{definition}\label{DefVeryRelativeBoundaries} For a map ${f : Y \rightarrow X}$ and subobjects $u$, $v$ of $X$, we write
\[ \partial_f u \leq v \]
 if $f$ factors through the subobject ${ u \vee (u \Rightarrow v)}$. 
\end{definition}

The notation is compatible with the ordering of subobjects in the following sense.

\begin{lemma}[`Transitivity'] \label{LemPseudoTransitivity} Let ${f : Y \rightarrow X}$ and let  ${v \leq w}$ be subobjects of  $X$. For every subobject $u$ of $X$,  ${\partial_f u \leq v}$ implies ${\partial_f u \leq w}$.
\end{lemma}
\begin{proof}
Since ${v \leq w}$, ${(u \Rightarrow v) \leq (u \Rightarrow w)}$ and so ${u \vee (u \Rightarrow v ) \leq u \vee (u \Rightarrow w)}$.
\end{proof}

\begin{definition}\label{DefSubobjecton top of Ell} A map $f$ with codomain $X$ is {\em on top of} a subobject ${v}$ of $X$  if ${\partial_f u \leq v}$  for every subobject $u$ of $X$.
\end{definition}

The intuition behind the terminology is that the dimension of the image of $f$ is `just above' the dimension of $v$.
At the present level of generality this is somewhat vague, but when $f$ is the identity and $v$ is the skeleton of a level then the intuition is fairly accurate.
See the comments after Definition~\ref{DefOnTopOfNaturalMono}.

For future reference we state some simple consequences of the definition.
First observe that the `transitivity' of Lemma~\ref{LemPseudoTransitivity} easily lifts to the following. 

\begin{lemma}\label{LemOnTopIsMonotone} If ${f : Y \rightarrow X}$ is on top of $v$ and ${v\leq w}$ then $f$ is on top of $w$.
\end{lemma}

Further useful properties are stated below.

\begin{lemma}\label{LemNearlyTrivial} For maps ${f : Y \rightarrow X}$,  ${g : Z \rightarrow Y}$  and $v$  a subobject of $X$, the following hold:
\begin{enumerate}
\item (Sieve-property/right-ideal) If $f$ is on top of  $v$ then so is ${f g}$.
\item If ${g}$ is a regular epimorphism then,  ${f g}$ on top of  $v$ implies $f$ on top of  $v$.
\item The map $f$ is on top of $v$ if and only if the image of $f$ is on top of $v$.
\item\label{ItmForLeftIdeal} If $g$ is on top of  ${f^* v}$ then ${f g}$ is on top of  $v$.
\item If $f$ is monic  then:  ${g}$ is on top of  ${f^* v = f \cap v}$  if and only if  ${f g}$ is on top of  $v$.
\end{enumerate}
\end{lemma}
\begin{proof}
The first item is  trivial. 
The second follows from orthogonality.
The third follows from the previous two.
To prove the fourth let $u$ be a subobject of $X$.
By hypothesis, $g$ factors through 
\[ f^* u \vee (f^* u \Rightarrow f^* v) = f^* (u \vee (u \Rightarrow v)) \]
so ${f g}$ factors through ${u \vee (u \Rightarrow v)}$

Consider now  the final item.
One implication does not need $f$ monic and follows from the fourth item.
To prove the converse let ${w}$ be a subobject of $Y$.  By hypothesis ${fg}$ factors through ${\exists_f w  \vee (\exists_f w \Rightarrow v)}$ so
$g$ factors through 
\[ f^* [\exists_f w  \vee (\exists_f w \Rightarrow v)] = f^* (\exists_f w) \vee  (f^* (\exists_f w) \Rightarrow f^* v) \]
 and, since $f$ is monic, ${f^* (\exists_f w) =w}$. (See \cite[Proposition~{IV.6.3}]{maclane2} and the comments after the proof there.) 
 Hence, $g$ factors through ${w \vee (w\Rightarrow f^* v)}$.
\end{proof}

For the discussion below, recall that a subcategory is {\em mono-coreflective} if it is coreflective and the counit of the coreflection is monic.

Assume now that our fixed Heyting category is  equipped with a mono-coreflective subcategory with (monic) counit $\ell$.
We emphasize the following refinement of the terminology.

\begin{definition}\label{DefOnTopOfNaturalMono}
A map $f$  is said to be {\em on top of $\ell$} if $f$ is on top of the monic $\ell_X$ where $X$ is the codomain of $f$. 
\end{definition}

 If $l$ is a level with monic skeleta ${\ell_X : l_! (l^* X) \rightarrow X}$ then, for any object $X$, $id_X$ is on top of ${\ell_X}$ if and only if $X$ has $l$-skeletal boundaries in the sense of \cite[Definition~{3.1}]{Menni2019a}.

\begin{example}[$-\infty$-skeletal boundaries and `Boolean' objects]\label{ExBooleanObjects}
Consider level $-\infty$ of a topos $\calE$. The $-\infty$-skeleton of an object $X$ in $\calE$ is just the initial subobject ${0 \rightarrow X}$.
For a subobject $u$ of $X$, ${\partial_{id} u \leq 0}$ if and only if $u$ is complemented \cite[Example~{3.3}]{Menni2019a}.
So $X$ has $-\infty$-skeletal boundaries if and only if the Heyting algebra of subobjects of $X$ is  Boolean. 
\end{example}

Intuitively, objects with $-\infty$-skeletal boundaries are `discrete'.
This has a rigorous formulation in the context of Axiomatic Cohesion. 
Recall that a geometric morphism ${p : \calE \rightarrow\calS}$ is {\em pre-cohesive} if it is local, hyperconnected, essential and the leftmost adjoint ${p_! : \calE \rightarrow \calS}$ preserves finite products.  Intuitively, $\calE$ is a topos `of spaces' over a topos $\calS$ `of sets' and the inverse image ${p^* : \calS \rightarrow \calE}$ is the full subcategory of `discrete' spaces.
Notice that ${p^* \dashv p_* \dashv p^! : \calS \rightarrow \calE}$ is a level of $\calE$. Intuitively, it is level~0.

\begin{proposition}\label{PropBooleanIffDiscrete} If $\calS$ is a Boolean topos and ${p : \calE \rightarrow \calS}$ is  pre-cohesive  then ${p^* : \calS \rightarrow \calE}$ coincides with the full subcategory of objects with  $-\infty$-skeletal boundaries.
\end{proposition}
\begin{proof}
The Aufhebung of level~${-\infty}$ (i.e. level~0) coincides with ${p^* \dashv p_* \dashv p^! : \calS \rightarrow \calE}$ by \cite[Corollary~{4.15}]{LawvereMenni2015}. The  full subcategory of objects with  $-\infty$-skeletal boundaries coincides with that of Boolean objects by Example~\ref{ExBooleanObjects} and so, it is included in ${p^* : \calS \rightarrow \calE}$ by \cite[Proposition~{6.7}]{MarmolejoMenni2021}.
Finally, since $\calS$ is Boolean and $p$ is hyperconnected, ${p^* A}$ is Boolean for every $A$ in $\calS$,
so ${p^* }$ factors through the subcategory of Boolean objects.
\end{proof}

In other words, an object is $0$-skeletal if and only if it has $-\infty$-skeletal boundaries.
We stress that, as observed in \cite[Example~{1.5}]{Menni2019a}, objects with 0-discrete boundaries need not be 1-skeletal.
Instead, in many examples, an object is 1-skeletal if and only if it is a quotient of a 0-separated object with discrete boundaries \cite[Proposition~{7.3}]{Menni2019a}.

We now continue studying the right ideal of maps on top of $\ell$.

\begin{lemma}\label{LemLeftIdeal} 
For every  ${f : Y \rightarrow X}$ and ${g : Z \rightarrow Y}$, if $g$ is on top of $\ell$ then so is ${f g}$.
\end{lemma}
\begin{proof}
Let $w$ be a subobject of $Y$.
Then ${w \vee (w \Rightarrow \ell_Y) \leq w \vee (w\Rightarrow f^* \ell_X)}$ as subobjects of $Y$  because ${\ell_Y \leq f^* \ell_X}$ by naturality of $\ell$.
Since  $g$ is on top of $\ell_Y$, we  have ${\partial_g w \leq \ell_Y \leq f^* \ell_X}$ and so, by Lemma~\ref{LemOnTopIsMonotone}, 
${\partial_g w \leq f^* \ell_X}$. Hence, $g$ is on top of ${f^* \ell_X}$ and Lemma~\ref{LemNearlyTrivial}(\ref{ItmForLeftIdeal}) implies that ${f g}$ is on top of $\ell_X$.
\end{proof}

\begin{proposition}\label{PropIdealOfMapsInTopOfell}
The maps on top of $\ell$ form an ideal. 
If, moreover, $\ell$ is mono-cartesian then the ideal is idempotent; in fact, for any $f$ on top of $\ell$ both maps in the epi/mono factorization of $f$ are on top of $\ell$.
\end{proposition}
\begin{proof}
The collection is a right-ideal by Lemma~\ref{LemNearlyTrivial} and it is a left-ideal by Lemma~\ref{LemLeftIdeal}.
To prove idempotency,  assume that $\ell$ is mono-cartesian and let ${f : Y  \rightarrow X}$ be on top of $\ell$.
Since the underlying category is Heyting by hypothesis (and hence regular), there are regular-epi/mono factorizations.
So we can let ${f = m e}$ with $m$ monic and $e$ a regular-epimorphism.
Then the image $m$ of $f$ is on top of $\ell$ by Lemma~\ref{LemNearlyTrivial}.
Also, as $m$ is monic and ${f = m e}$ is on top of $\ell$, $e$ is on top of ${m^* \ell_X}$ by Lemma~\ref{LemNearlyTrivial}.
Since $\ell$ is mono-cartesian, $e$ is on top of ${\ell_{M} = m^* \ell_X }$ where $M$ is the domain of $m$. 
\end{proof}

Roughly speaking, given an idempotent comonad with monic and mono-cartesian counit on a Heyting category  we have `derived' an idempotent ideal in the same category. If the underlying Heyting category is a presheaf topos then, under certain conditions, we will be able to restrict the derived ideal to an idempotent ideal in the site.

To give more details let $\calD$ be a small category.
If an idempotent comonad on ${\Psh{\calD}}$ has a monic counit $\ell$ then we say that  a map $f$ in $\calD$ is {\em on top of $\ell$} if ${\calD(-,f)}$ is on top of $\ell$ in the sense of Definition~\ref{DefOnTopOfNaturalMono}.
This is a natural extension, in the present context, of the terminology we have been using.

\begin{lemma}\label{LemRestrictedIdealOfMapsOnTopOfNaturalMono} 
If $\calD$ is a small category with split-epi/mono factorizations then, for every  mono-coreflective subcategory of ${\Psh{\calD}}$ with mono-cartesian  counit $\ell$, the maps  in $\calD$ that are on top of $\ell$  form a mono-cartesian idempotent ideal in $\calD$.
\end{lemma}
 \begin{proof}
 The maps in ${\Psh{\calD}}$ that are on top of $\ell$ form an idempotent ideal by Proposition~\ref{PropIdealOfMapsInTopOfell}.
  It easily follows that the  maps $f$ in $\calD$  on top of $\ell$  form an ideal in $\calD$.
To prove that the latter ideal is idempotent requires a bit more.
Assume that ${f : E \rightarrow D}$ in $\calD$ is on top of $\ell$. 
Let ${f = m e}$ with $m$ monic and $e$ in split epic.
Then  ${\calD(-,f) = \calD(-,m) \calD(-,e)}$ is the usual factorization in $\Psh{\calD}$ as a monic map following an epic one. 
So, by  Proposition~\ref{PropIdealOfMapsInTopOfell} again, both ${\calD(-,m)}$ and ${\calD(-,e)}$ are on top of $\ell$.
Hence $e$ and $m$ are on top of $\ell$. This completes the proof that the ideal is idempotent but, moreover, the map $e$ has a section by hypothesis,  so it also follows from  Lemma~\ref{LemMonoCartesian} that the ideal is mono-cartesian.
 \end{proof}

Roughly speaking, for $\calD$ as above, we may restrict mono-cartesian idempotent ideals along Yoneda.
In particular, we may apply Lemma~\ref{LemRestrictedIdealOfMapsOnTopOfNaturalMono} to the principal comonads induced by mono-cartesian ideals which, as we already know by Proposition~\ref{PropMonoCartesianIdeals} are in correspondence with the  full  subcategories of $\calD$ that are closed under subobjects.

\begin{remark}[The successor of a subcategory closed under subobjects]\label{RemSuccessorCat} Assume that  $\calD$ has split-epi/mono factorizations and let ${i : \calC \rightarrow \calD}$  a full  subcategory closed under subobjects.  Proposition~\ref{PropMonoCartesianIdeals} implies that the corresponding ideal $\calI$ in $\calD$ is mono-cartesian. That is, the counit $\ell$ of the induced principal comonad on $\Psh{\calD}$ is mono-cartesian.  Let $\calJ$ be the idempotent ideal in ${\Psh{\calD}}$ consisting of maps on top of $\ell$.
By Lemma~\ref{LemRestrictedIdealOfMapsOnTopOfNaturalMono}, $\calJ$ restricts to a mono-cartesian ideal $\calI'$ in $\calD$.
Again by Proposition~\ref{PropMonoCartesianIdeals}, $\calI'$ is uniquely determined by the full subcateory ${i' : \calC' \rightarrow \calD}$ consisting of the objects $D$ in $\calD$ such that $id_D$ is in $\calJ$, that is, such that $id_D$ is on top of $\ell$. This subcategory $i'$ will be called the {\em successor} of $i$.
\end{remark}

The successor of a subcategory as above deserves a description purely in terms of $\calD$.
To give one we first need more specific information about maps on top of principal comonad.

\section{Maps on top of a principal comonad}
\label{SecMapsOnTopOfPrincipalComonad}

Let $\calD$ be a small category, let $\calI$ be an idempotent ideal therein,  and let ${\ell}$ be the (monic) counit of the  induced  principal comonad on ${\Psh{\calD}}$.

For any subobject ${u : U \rightarrow X}$ in $\Psh{\calD}$, the subobject ${u \Rightarrow \ell : (U \Rightarrow \ell) \rightarrow X}$ may be described as follows
\[ (U \Rightarrow \ell) D = \{ x \in X D \mid \textnormal{for all } f : E \rightarrow D  , \  x \cdot f \in U E \textnormal{ implies  } x \cdot f \textnormal{ is $\calI$-generated} \} \subseteq X D \]
for each object $D$ in $\calD$.

The next result is both and abstraction and a refinement of \cite[Lemma~{5.1}]{Menni2019a}.
In order to state it we need to introduce some notation.
For each $K$ in $\calD$ and ${k \in X K}$, the image of the corresponding  ${k : \calD(-, K) \rightarrow  X}$ will be denoted by ${u_k : U_k \rightarrow X}$ so that, for every $J$ in $\calD$ and ${j \in X J}$, ${j \in U_k J \subseteq X J}$ if and only if there is a ${g : J \rightarrow K}$ such that ${k \cdot g = j}$.

\begin{lemma}\label{LemAnalogueOf5.1}  
For any $X$ in $\Psh{\calD}$, $D$ in $\calD$ and ${x \in X D}$, the following are equivalent:
\begin{enumerate}
\item The map ${x : \calD(-, D) \rightarrow X}$  is on top of  $\ell_X$.
\item  For every ${f : E \rightarrow D}$ in $\calD$, ${\partial_x u_{x \cdot f} \leq \ell_X}$.
\item For every ${f : E \rightarrow D}$, either ${x \cdot f}$ is $\calI$-generated or there is a ${g : D \rightarrow E}$ in $\calD$ such that ${x \cdot (f g) = x}$.
\end{enumerate}
\end{lemma}
\begin{proof}
The first item trivially implies the second.
To prove that the second implies the third let ${f : E \rightarrow D}$ in $\calD$ and consider the subobject  ${u_{x \cdot f} : U_{x \cdot f} \rightarrow X}$.
By hypothesis, ${x \in U_{x\cdot f} D}$ or ${x \in (U_{x \cdot f} \Rightarrow \ell)  D}$.
If ${x \in  U_{x\cdot f} D}$ then there is a ${g : D \rightarrow E}$ such that ${(x \cdot f)\cdot g = x}$.
On the hand, if ${x \in (U_{x\cdot f} \Rightarrow \ell)  D}$ then, for all ${h : B \rightarrow D}$ in $\calD$, ${x \cdot h \in U_{x\cdot f} B}$ implies ${x \cdot h}$ is $\calI$-generated. Taking ${h = f :  E \rightarrow D}$ we  conclude that ${x\cdot f}$ is $\calI$-generated.

To prove that the third item implies the first, let ${u : U \rightarrow X}$ be a subobject and assume that ${x \not\in U D \subseteq X D}$.
We show that ${x \in (U \Rightarrow \ell) D \subseteq X D}$.
To do this let ${f : E \rightarrow D}$ be such that ${x \cdot f \in U E}$.
Then, by hypothesis, either ${x \cdot f}$ is $\calI$-generated  or there is a ${g : D \rightarrow E}$ in $\calD$ such that ${x \cdot (f g) = x}$.
If such a $g$ exists then, as ${x \cdot f \in U E}$, ${(x \cdot f)\cdot g \in U D}$. So ${x \in U D}$, which contradicts our assumption.
Hence, ${x \cdot f \in U E}$ implies that ${x \cdot f}$ is $\calI$-generated.
That is, ${x \in (U \Rightarrow \ell) D \subseteq X D}$ as we needed to show.
\end{proof}

Specializing to representables we obtain the following analogue of  \cite[Lemma~{5.4}]{Menni2019a}.

\begin{lemma}\label{LemPreAnalogueOf5.4} 
For every map ${e : E \rightarrow D}$ in $\calD$ the following hold:
\begin{enumerate}
\item $e$ is on top of $\ell$   if and only if, for every ${f : F \rightarrow E}$ in $\calD$, either ${e f}$ is in $\calI$ or there is a ${g :E \rightarrow F}$ such that ${e f g = e}$.
\item If $e$ is monic then, $e$ is on top of $\ell$ if and only if, for every ${f : F \rightarrow E}$ in $\calD$, either ${e f}$ is in $\calI$ or $f$ has a section.
\item If $e$ has a retraction then, $e$ is on top of $\ell$ if and only if, for every $f$ with codomain $E$, $f$ is in $\calI$ or $f$ has a section.
\end{enumerate}
\end{lemma}
\begin{proof}
By  Lemma~\ref{LemAnalogueOf5.1}, $d$ is on top of $\ell$ if and only if for every ${f : F \rightarrow E}$ in $\calD$, either ${e f}$ is $\calI$-generated  or there is a ${g :E \rightarrow F}$ such that ${e f g = e}$.
By the left-ideal property, ${e f}$ is $\calI$-generated if and only if ${e f}$ is in $\calI$.

If $e$ is monic,  ${e f g = e}$ is equivalent to ${f g = id}$, so the second item holds.
Finally, assume that $e$ has a retraction.
It remains to show that  ${e f}$ in $\calI$ implies $f$ in $\calI$,
but this follows from existence of a retraction for $e$ and the left ideal property.
\end{proof}

We may now return to the issue raised in Remark~\ref{RemSuccessorCat}.

\begin{theorem}\label{ThmSuccessorCats} 
Let $\calD$ be a  small category with split-epi/mono factorizations.
If ${\calC \rightarrow \calD}$ is a full  subcategory closed under subobjects then its successor is the subcategory of $\calD$ consisting of the objects $D$ in $\calD$ such that:  for every monic ${f : C \rightarrow D}$ with codomain $D$, $C$ is in  $\calC$ or $f$ is an isomorphism.
In particular, as a function on the poset of (full and closed-under-subobjects) subcategories of $\calD$, the successor is monotone increasing.
\end{theorem}
\begin{proof}
As already observed in Remark~\ref{RemSuccessorCat}, the successor consists of the objects $D$ in $\calD$ such that $id_D$ is in $\calJ$, that is, such that $id_D$ is on top of $\ell$. By Lemma~\ref{LemPreAnalogueOf5.4}, this holds if and only if, for every $f$ in $\calD$ with codomain $D$, $f$ factors through $\calC$ or $f$ has a section. This is equivalent, under the assumption of split-epic/mono factorizations, to the condition in the statement.
\end{proof}

\section{Examples}
\label{SecExamples}

In this section we let $\calD$ be a small category with split-epi/mono factorizations.

\begin{example}\label{ExPosets} 
Let $\calD$ be poset. The full subcategories that are closed under subobjects are just the lower-sets in $\calD$. 
The successor of a lower set ${\calC \rightarrow \calD}$ consists of those $D$ in $\calD$ such that, for every ${C \leq D}$, $C$ is in $\calC$ or ${C = D}$.   
\end{example}

In particular, the successor of the empty lower-set  in a poset  is the lower-set of minimal elements, and the successor of the largest ideal is itself.
Non-surprisingly,  this is also a particular case of the following more general observation.

\begin{example}[Extremes]\label{ExSuccesorOfEmpty}
The successor of the empty subcategory of $\calD$ is the subcategory determined by the objects $D$ such that every monomorphism with codomain $D$ is an isomorphism. 
The successor of the identity ${\calD \rightarrow \calD}$ is itself. 
\end{example}

\begin{example}[The classifier of non-trivial Boolean algebras]\label{ExBAclassifier}
Let $\mathbb{F}$ be the category of non-empty finite sets.
For each $N$ in $\mathbb{F}$ let ${\mathbb{F}_N \rightarrow \mathbb{F}}$ be the full subcategory of those $M$ such that there is a monic ${M \rightarrow N}$. The full subcategories of $\mathbb{F}$ that are closed under subobjects are  the ${\mathbb{F}_N \rightarrow \mathbb{F}}$, together with the empty subcategory and the whole subcategory.
The successor of $\emptyset$ is ${\mathbb{F}_1}$ by Example~\ref{ExSuccesorOfEmpty}.
We next prove that the succesor of ${\mathbb{F}_N}$ is ${\mathbb{F}_{N+1}}$.
Every subobject of  ${N + 1}$ has domain in  ${\mathbb{F}_N}$ or it is an iso.
That is, ${N + 1}$ is in the successor.
On the other hand, the obvious coproduct injection ${N + 1 \rightarrow N + 1 +1 = N+ 2 }$ is not an iso and its domain is not in ${\mathbb{F}_N}$.
Hence, ${N + 2}$ is not in the successor of ${\mathbb{F}_N}$.
Any subcategory of $\mathbb{F}$ that is closed under subobjects and properly containing ${\mathbb{F}_N}$ is either ${\mathbb{F}_{N + 1}}$ or it contains ${N + 2}$. So the  successor of ${\mathbb{F}_N}$ is ${\mathbb{F}_{N + 1}}$.
\end{example}

The calculation of the Aufhebung for ${\Psh{\mathbb{F}}}$ remains an open problem \cite{LawvereOpenProblems}.
The next example, together with \cite[Theorem~{3.21}]{KennetEtAl2011}, shows that the successor is  different from the Aufhebung in general.

\begin{example}[Simplicial sets]\label{ExSimplicialSets}
Let ${\Delta}$ be the category of finite non-empty totally ordered sets.
As usual, for any ${d \in \Nat}$,  we let ${[d]}$ be the object in $\Delta$ with ${d + 1}$  elements.
Also, we let ${\Delta_d \rightarrow \Delta}$ be the full subcategory consisting of the objects ${[c]}$ with ${c\leq d}$.
Just as in Example~\ref{ExBAclassifier}, the subcategories of $\Delta$ closed under subobjects are the extreme ones together with those of the form ${\Delta_d \rightarrow \Delta}$ for some ${d \in \Nat}$. 
Moreover, essentially the same argument in Example~\ref{ExBAclassifier} implies that the successor of ${\Delta_d}$ is ${\Delta_{d+1}}$.
\end{example}

We are mainly interested in toposes of spaces. 
It is therefore natural to sometimes assume that ${\Psh{\calD}}$ is pre-cohesive \cite{LawvereMenni2015}.
So consider the following variant of Proposition~\ref{PropBooleanIffDiscrete}.

\begin{proposition}\label{PropAufEqualsSuc} 
If $\calD$ has a terminal object and every object has a point (so that ${\Psh{\calD}}$ is a pre-cohesive over $\Set$), then the successor of level ${-\infty}$ coincides with its Aufhebung, namely, the subtopos of $\neg\neg$-sheaves (which is a level in this case).
\end{proposition}
\begin{proof}
The Aufhebung of Level~${-\infty}$  of ${\Psh{\calD}}$ coincides with the subtopos of $\neg\neg$-sheaves by \cite[Corollary~{4.5}]{LawvereMenni2015}. Under the present hypotheses, the $\neg\neg$-topology (as described in \cite[II.2(d)]{maclane2}) is such that a sieve on $D$ covers if and only if it contains all the points of $D$. This is obviously a rigid topology and the corresponding subtopos is that induced by the full subcategory consisting of the terminal object. This subcategory is,  since every object of $\calD$ has a point, the succesor of the empty subcategory by Example~\ref{ExSuccesorOfEmpty}. 
\end{proof}

Let ${\calD_0 \rightarrow \calD}$ be the successor of the empty subcategory (Example~\ref{ExSuccesorOfEmpty}) and  let ${\calD_1 \rightarrow \calD}$ be the successor of ${\calD_0 \rightarrow \calD}$. By Theorem~\ref{ThmSuccessorCats} and Example~\ref{ExSuccesorOfEmpty} we infer:

\begin{lemma}\label{LemSucOfzero} An object $D$ in $\calD$ is in ${\calD_1 \rightarrow \calD}$ if and only if for every monic ${f : C \rightarrow D}$ in $\calD$, $f$ is an isomorphism or every monic with codomain $C$ is an isomorphism.
\end{lemma}

Even if ${\Psh{\calD}}$ is pre-cohesive, the successor of 0 need not coincide with its Aufhebung.

\begin{example}[Level~${\epsilon}$ cannot be the Aufhebung of $0$ but it may be the successor of $0$]
Assume that $\calD$ is finite, that it has a terminal object and that every object has a point. 
By \cite[Corollaries~{5.4} and {4.5}]{MMlevelEpsilon}, the  full subcatgory ${\calD_! \rightarrow \calD}$   of objects that have exactly one point determines a level that is not way-above level~0, so it cannot be the Aufhebung of $0$.
In particular, consider the 4-element graphic monoid discussed in \cite[p.~62]{Lawvere1989taco} and let ${\calD}$ be its idempotent splitting.
As explained in \cite{MMlevelEpsilon} the idempotent splitting looks as follows
\[\xymatrix{
 & D \ar[d]<-1ex>_-s  \\
1 \ar@(u,l)[ru]^-{\ddagger} \ar[r]<-1ex>_-{\bot}  \ar[r]<+1ex>^-{\top} &  G \ar@(ru,rd)^-{\alpha} \ar[u]<-1ex>_-r
}\]
(without drawing the maps towards 1 and) with ${\alpha \alpha = \alpha = s r}$, ${r s = id}$, ${r \bot = \ddagger  = r \top}$ and ${s \ddagger = \bot}$.
In this case, the subcategory ${\calD_! \rightarrow \calD}$ coincides with ${\calD_1 \rightarrow \calD}$, because $D$ has exactly two subobjects: $\ddagger$ and the identity. So the level ${\Psh{\calD_!} \rightarrow \Psh{\calD}}$ is the successor of level~0 but, as  mentioned above, it cannot be the Aufhebung of level~0. In this example the poset of full subcategories of $\calD$ closed under subobjects looks as follows
\[\xymatrix{
-\infty  = \emptyset < \calD_0 = \{ 1 \} < \calD_1 = \{ 1, D \} < \calD_2 = \calD 
}\]
and the successor is the obvious one.
\end{example}

In some cases  the succesor and the Aufhebung of 0  coincide.
Consistently with \cite{Menni2019a}, where objects in $\calD_1$ are called {\em edge types}, we say that an object in $\calD$ is {\em edge-wise connected} if for every parallel pair of points ${u, v : 1 \rightarrow D}$ there exists an ${E}$ in $\calD_1$ and a map ${f : E \rightarrow D}$ such that both $u$ and $v$ factor through $f$.

\begin{corollary}\label{CorSuc0}
Assume that  $\calD$  has a terminal object and that  every object has a point (so that ${\Psh{\calD}}$ is pre-cohesive).
If $1$ is the only object in ${\calD_1}$ with exactly one point and every object in $\calD$ is edge-wise connected then the Aufhebung of level 0 coincides with ${\Psh{\calD_1} \rightarrow \Psh{\calD}}$, the successor of level 0.
\end{corollary}
\begin{proof}
Follows from \cite[Proposition~{7.2}]{Menni2019a} which, under slightly more general hypotheses (no need of split-epi/mono factorizations) shows that the Aufhebung of level 0 is ${\Psh{\calD_1} \rightarrow \Psh{\calD}}$ which, in this paper, we have identified as the successor of $0$.
\end{proof}

If we allow ourselves to picture an object with a unique point as `a point with infinitesimal information around it' then Corollary~\ref{CorSuc0} says, intuitively, that if there `no infinitesimals' and the objects of ${\calD_1}$ are enough to connect points then the Aufhebung of 0 coincides with the successor of 0.

Further results relating the Aufhebung and the successor may help to illuminate the nature of both, and maybe aid in the calculation of the Aufhebung in some cases. We should also mention that Lawvere defines in \cite{Lawvere91b} an {\em addition} of levels and invites to compare the Aufhebung with the function that results from adding level~1. This addition-with-level-1 should also be compared with the successor presented here.

From now on let $\calE$ be another category with  split-epic/mono factorizations so that  ${\calD \times \calE}$ also has split-epic/mono factorizations.

\begin{lemma}\label{LemSucMinusInftyInProduct} The subcategory  ${(\calD \times \calE)_0 \rightarrow \calD \times \calE}$ coincides with  ${\calD_0 \times \calE_0 \rightarrow \calD\times \calE}$
\end{lemma}
\begin{proof}
Every monic with codomain ${(D, E)}$ in ${\calD \times \calE}$ is an iso if and only if the same holds for $D$ and for $E$.
So the result follows again from  Example~\ref{ExSuccesorOfEmpty}.
\end{proof}

As expected, products with the successor of the empty category have special properties.

\begin{lemma}\label{LemLeibnizPropertyOfSuccessor} 
With $\calD$ and $\calE$ as above,  let ${\calC \rightarrow \calD}$ a full subcategory closed under subobjects and let  ${\calC' \rightarrow \calD}$ be its successor. Then the successor of ${\calC \times \calE_0 \rightarrow \calD \times \calE}$ contains ${\calC' \times \calE_0 \rightarrow \calD \times \calE}$.
\end{lemma}
\begin{proof}
Let ${(D, E)}$ be in ${\calC' \times \calE_0 \rightarrow \calD \times \calE}$.
Let ${(a, b) : (A, B) \rightarrow (D,E)}$ be a monic in  ${\calD \times \calE}$.
That is ${a : A \rightarrow D}$ is monic in $\calD$ and ${b : B \rightarrow E}$ is monic in $\calE$.
Since $E$ is in $\calE_0$, $b$ is an isomorphism.
Since $D$ is in $\calC'$,  $A$ is in $\calC$ or $a$ is an isomorphism.
If $A$ is in $\calC$ then ${(A, B)}$ is in ${\calC\times \calE_0}$.
If $a$ is an isomorphism then so is ${(a, b)}$.
Altogether ${(D, E)}$ is in the successor of ${\calC \times \calE_0 \rightarrow \calD \times \calE}$.
\end{proof}

Hence, the Leibniz rule in the next result may not be a surprise.

\begin{proposition}\label{Prop1inProductCat} 
The subcategory  ${(\calD \times \calE)_1 \rightarrow \calD\times \calE}$ coincides with the union of   ${\calD_1 \times \calE_0 \rightarrow \calD \times \calE}$ and ${\calD_0 \times \calE_1 \rightarrow \calD \times \calE}$.
\end{proposition}
\begin{proof}
Lemmas~\ref{LemSucMinusInftyInProduct} and \ref{LemLeibnizPropertyOfSuccessor} imply that ${(\calD \times \calE)_1 \rightarrow \calD\times \calE}$  contains the indicated union.
Now let ${(D, E)}$ be an arbitrary object of ${(\calD \times \calE)_1}$.
Lemma~\ref{LemSucOfzero} implies that, for every monic ${(a, b) : (A, B) \rightarrow (D, E)}$, ${(a, b)}$ is an iso or ${(A, B)}$ is in ${(\calD \times \calE)_0}$.  Equivalently, by Lemma~\ref{LemSucMinusInftyInProduct}, either both $a$, $b$ are isos or:  both $A$ is in $\calD_0$ and $B$ is in $\calE_0$.

If neither $D$ nor $E$ have a proper subobject then ${(D, E) \in \calD_0 \times \calE_0}$  by Example~\ref{ExSuccesorOfEmpty}, so ${(D, E)}$  is  both in ${\calD_1 \times \calE_0}$ and ${\calD_0 \times \calE_1}$.
Otherwise there is a proper subobject of $D$ or a proper subobject of $E$.
If there is a non-iso monic ${a : A \rightarrow D}$ in $\calD$ then ${(a, id_E) : (A, E) \rightarrow (D, E)}$ is a non-iso monic.
Hence, both $A$ in $\calD_0$ and $E$ in $\calE_0$.
Hence, if there is a non-iso monic ${A \rightarrow D}$, then $D$ is in $\calD_1$ and $E$ is in $\calE_0$.
Similarly, if there is a non-iso monic ${B \rightarrow E}$ then $E$ is in $\calE_1$ and $D$ is in $\calD_0$. 
\end{proof}

Let ${(\calD \times \calE)_2 \rightarrow \calD \times \calE}$ be the successor of  ${(\calD \times \calE)_1 \rightarrow \calD \times \calE}$.

\begin{proposition}\label{LemFor2} 
The subcategory   ${(\calD \times \calE)_2 \rightarrow \calD \times \calE}$ coincides with ${\calD_1 \times \calE_1 \rightarrow \calD \times \calE}$.
\end{proposition}
\begin{proof}
Let ${(D, E)}$ be in  ${\calD_1 \times \calE_1}$ and let ${(m, n) : (A, B) \rightarrow (D, E)}$ be monic.
Equivalently, both ${m : A \rightarrow D}$ and ${n : B \rightarrow E}$ are monic. 
So the following two items hold:
\begin{enumerate}
\item $m$ is an iso or $A$ is in $\calD_0$,
\item $n$ is an iso or $B$ is in ${\calE_0}$.
\end{enumerate}
Equivalently, one of the following items holds:
\begin{enumerate}
\item both $m$ and $n$ are isos or
\item $m$ is iso and $B$ is in $\calE_0$ or
\item $A$ is in $\calD_0$ and $n$ is an iso or
\item $A$ is in $\calD_0$ and $B$ is in $\calE_0$.
\end{enumerate}
So, ${(m,n)}$ is an iso, or ${(A, B)}$ is in ${\calD_1 \times \calE_0}$, or ${(A, B)}$ is in ${\calD_0 \times \calE_1}$.
Equivalently,  by Proposition~\ref{Prop1inProductCat}, ${(m,n)}$ is an iso, or ${(A, B)}$ is in ${(\calD \times \calE)_1 }$.
Therefore, ${(D, E)}$ is in ${(\calD \times \calE)_2 \rightarrow \calD\times \calE}$.

On the other hand, let ${(D, E)}$ be in ${(\calD \times \calE)_2}$.
That is, for every monomorphism  ${(m, n) : (A, B) \rightarrow (D, E)}$ in ${\calD\times \calE}$, ${(m, n)}$ is an iso or  ${(A, B)}$ is in ${(\calD \times \calE)_1}$.
Equivalently, by  Proposition~\ref{Prop1inProductCat}, ${(m,n)}$ is an iso, or ${(A, B)}$ is in ${\calD_1 \times \calE_0}$, or ${(A, B)}$ is in ${\calD_0 \times \calE_1}$.
We next prove that $D$ is in $\calD_1$.
If every monic with codomain $D$ is an iso then $D$ is in $\calD_0$ and, therefore, in $\calD_1$.
So let ${m : A \rightarrow D}$ be a non-iso monic. Then ${(m, id) : (A, E) \rightarrow (D, E)}$ is not an iso.
Hence, ${(A, E)}$ is in ${\calD_1 \times \calE_0}$, or ${(A, E)}$ is in ${\calD_0 \times \calE_1}$. 
In any case, $A$ is in $\calD_1$.
Similarly, $E$ is in $\calE_1$.
\end{proof}

The idea does not persist  above ${(\calD \times \calE)_2 \rightarrow \calD \times \calE}$ in an obvious way.
Consider  ${\mathbb{F}\times \mathbb{F}}$ as an example.
It follows from the results above that  ${(\mathbb{F}\times \mathbb{F})_0 = \mathbb{F}_0 \times \mathbb{F}_0}$, that
 ${(\mathbb{F}\times \mathbb{F})_1}$ is the union of the subcategories ${\mathbb{F}_1 \times \mathbb{F}_0}$ and ${\mathbb{F}_0 \times \mathbb{F}_1}$, and that ${(\mathbb{F}\times \mathbb{F})_2 = \mathbb{F}_1 \times \mathbb{F}_1}$.

\begin{proposition}\label{Prop2and3} The successor of ${(\mathbb{F}\times \mathbb{F})_2 \rightarrow \mathbb{F}\times \mathbb{F}}$ is itself.
\end{proposition}
\begin{proof}
We prove something more general. For any ${(M, N)}$ in ${ \mathbb{F}\times \mathbb{F}}$, let  ${\calC_{M, N} \rightarrow  \mathbb{F}^2}$ be the full subcategory consisting of the objects that appear as the domain of a monic into ${(M, N)}$.
We give a sufficient condition  for the successor of ${\calC_{M, N}}$ to be itself.
Let ${(M, N)}$ be such that there are monos ${ m : 2 \rightarrow M}$ and ${n : 2 \rightarrow N}$.
  Consider the object ${(M +1, N)}$. It is easy to produce a monic  ${(id, y) : (M + 1, 1) \rightarrow (M+1, N)}$ which is not an iso because ${y  : 1 \rightarrow N}$ is not. Also,  ${(M+1, 1)}$ is not in ${\calC_{M, N}}$. So ${(M +1, N)}$ is not in the successor of ${\calC_{M, N} \rightarrow  \mathbb{F}^2}$. Similarly ${(M, N+1)}$ is not in that successor either.
 Moreover any object not in ${\calC_{M, N} \rightarrow  \mathbb{F}^2}$  is the codomain of an monic with domain ${(M+1, N)}$ or ${(M, N+1)}$.
 So the successor of ${\calC_{M, N} \rightarrow  \mathbb{F}^2}$ is  itself.
\end{proof}

The concept of object with skeletal boundaries seems to have wider significance. 
As a minor additional remark we calculate the objects on top of a conspicuous object in a small category without split-epi/mono factorizations.

\begin{example}
Consider the category of finite trees  as a full subcategory of that of posets.
(There is an obvious surjection 
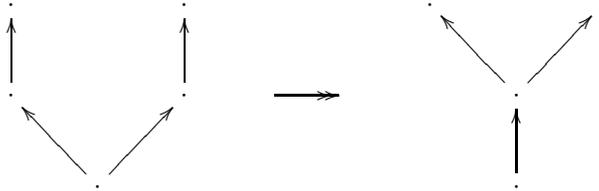
\[\xymatrix{
\cdot & & \cdot                      &                        &&  \cdot & &  \cdot \\
\cdot \ar[u] & & \cdot \ar[u] & \ar@{->>}[r] &&  & \ar[lu] \cdot \ar[ru] \\
 & \ar[lu] \cdot \ar[ru]           &                        &&  & & \cdot \ar[u]
}\]
which does not have a section.)
Consider the ideal of maps that factor through the tree  ${\uparrow}$ with two nodes,
and also the  trees whose identities are on top of (the counit of the principal comonad determined by) that  ideal.
Using Lemma~\ref{LemPreAnalogueOf5.4} it is straightforward to check that the identities of the following  two trees
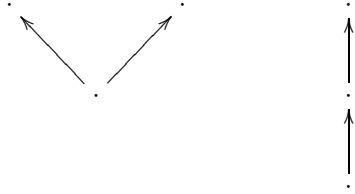
\[\xymatrix{
\cdot &                                 & \cdot &&  \cdot \\
         & \ar[lu] \cdot \ar[ru] &           && \cdot \ar[u] \\
         &                                 &           && \cdot \ar[u] 
}\]
are indeed on top of the ideal. Any other tree with more nodes would contain one of the above as a subtree, so the corresponding identity could not be on top of the ideal.
\end{example}

\end{document}